\documentclass[12pt,a4paper]{article}
\usepackage{amsmath,amsfonts,amssymb,amscd}
\usepackage[body={6.6in, 9.2in},left=1in, top=1.2in]{geometry}

\def\id{\operatorname{id}}
\def\ev{\operatorname{ev}}

\def\d{\operatorname{d}}

\def\pr{\operatorname{pr}}
\def\Im{\operatorname{Im}}
\def\pt{\operatorname{pt}}

\newcounter{th}
\def\t{\refstepcounter{th}{\bf \noindent{Theorem} \arabic{th}. }}

\newcounter{le}
\def\l{\refstepcounter{le}{\bf \noindent{Lemma} \arabic{le}. }}

\begin{document}


\begin{center}
\large{\bf ON THE STRUCTURE OF COMPLEX HOMOGENEOUS SUPERMANIFOLDS }
\end{center}

\bigskip

\begin{center}
    E.G. Vishnyakova \footnote{Mathematics Classification Number 32C11

    \textit{Key words and phrases.} Complex homogeneous supermanifold, Lie supergroup, action of a Lie
    supergroup.

 Work supported by German Academic Exchange Service, Michail-Lomonosov-Forschungsstipendien and by the Russian Foundation for
Basic Research (grant no. 07-01-00230).
 }

\bigskip
Tver State University

\end{center}

\bigskip

\noindent\textsc{Abstract.} For a Lie group $G$ and a closed Lie
subgroup $H\subset G$, it is well known that the coset space $G/H$
can be equipped with the structure of a manifold homogeneous under
$G$ and that any $G$-homogeneous manifold is isomorphic to one of
this kind. An interesting problem is to find an analogue of this
result in the case of supermanifolds.

In the classical setting, $G$ is a real or a complex Lie group and
$G/H$ is a real and, respectively, a complex manifold. Now, if $G$
is a real Lie supergroup and $H\subset G$ is a closed Lie
subsupergroup, there is a natural way to consider $G/H$ as a
supermanfold. Furthermore, any $G$-homogeneous real supermanifold
can be obtained in this way, see \cite{Kostant}. The goal of this
paper is to give a proof of this result in the complex case.

\bigskip

\bigskip

\begin{center}
{\bf 1. Preliminaries}
\end{center}

\bigskip

We will use the word "supermanifold" \,in the sense of
Beresin-Le\v{i}tes-Kostant (see \cite{Kostant, ley}). All the time,
we will be interested in the complex-analytic version of the theory.
A morphism $(M,\mathcal O_M)\to (N,\mathcal O_N)$ between two
supermanifolds is denoted by $\varphi=(\varphi_1,\varphi_2)$, where
$\varphi_1:M\to N$ is a continuous mapping and $\varphi_2:\mathcal
O_N\to (\varphi_1)_*(\mathcal O_M)$ is a homomorphism of sheaves of
superalgebras. We begin with the more general notion of a Lie
supergroup.

A {\it Lie supergroup} is a supermanifold $(G,\mathcal O_G)$, for
which the following three morphisms are defined:
  $\nu:(G,\mathcal O_G)\times (G,\mathcal O_G)\to
(G,\mathcal O_G)$ (multiplication morphism), $\iota:(G,\mathcal
O_G)\to (G,\mathcal O_G)$ (passing to the inverse),
$\varepsilon:(\pt,\mathbb C)\to (G,\mathcal O_G)$ (identity
morphism). Moreover, these morphisms should satisfy the usual
conditions, modeling the group axioms. The underlying manifold $G$
is a Lie group. We denote by $\mathfrak{g}$ the Lie superalgebra of
$(G,\mathcal O_G)$ (see \cite{ley} for the corresponding
definition).

An {\it action of a Lie supergroup $(G,\mathcal O_G)$ on a
supermanifold} $(M,\mathcal O_M)$ is a morphism $\mu:(G,\mathcal
O_G)\times (M,\mathcal O_M)\to (M,\mathcal O_M)$, such that the
following conditions hold:
\begin{itemize}
  \item $\mu \circ (\nu\times \id)=\mu\circ (\id\times \mu)$;
  \item $\mu\circ (\varepsilon\times \id)=\id$.
\end{itemize}

\noindent We denote by $\mathfrak v(M,\mathcal O_M)$ the Lie
superalgebra of holomorphic vector fields on $(M,\mathcal O_M)$. Let
$\mathfrak m_x$ be the maximal ideal of the local superalgebra
$(\mathcal O_M)_x$. The vector superspace $T_x(M,\mathcal
O_M)=(\mathfrak m_x/\mathfrak m_x^2)^*$ is called the {\it tangent
space} to $(M,\mathcal O_M)$ at $x\in M$. From the inclusions
$v(\mathfrak m_x)\subset (\mathcal O_M)_x$ and $v(\mathfrak
m_x^2)\subset \mathfrak m_x$, where $v\in \mathfrak v(M,\mathcal
O_M)$, it follows that there exists an even linear mapping
$\ev_x(v):\mathfrak m_x/\mathfrak m_x^2\to (\mathcal
O_M)_x/\mathfrak m_x\simeq \Bbb C$. In other words, $\ev_x(v)\in
T_x(M,\mathcal O_M)$, and so we obtain a map $\ev_x:\mathfrak
v(M,\mathcal O_M)\to T_x(M,\mathcal O_M)$.

Let $(U,\mathcal O_M)\subset (M,\mathcal O_M)$ be a superdomain with
even and, respectively, odd coordinates $(x_i)$ and $(\xi_j)$ and
let $f\in \mathcal O_M(U)$. We can write $f$ in the form
$$
f=f_0+\sum_if_i\xi_i+ \sum_{i\,j}f_{i\,j}\xi_i\xi_j+\ldots,
$$
where $f_{ij\ldots}$ are some holomorphic functions on $U$. For
$p\in U$ we will denote by $f(p)$ the value of $f_0$ at $p$.

Let $X\in T_x(M,\mathcal O_M)$. There is a neighborhood $(U,\mathcal
O_M)$ of the point $x$ and a vector field $v_X\in
\mathfrak{v}(U,\mathcal O_M)$ such that $\ev_x(v_X)=X$. We can
consider $X$ as a linear function on $(\mathcal O_M)_x$. Namely,
 $X(f_x):=(v_X(f_x))(x)$, where
$f_x\in (\mathcal O_M)_x$. 

Let $\mu=(\mu_1,\mu_2):(G,\mathcal O_G)\times (M,\mathcal O_M)\to
(M,\mathcal O_M)$ be an action. Then there is a homomorphism of the
Lie superalgebras $\widetilde{\mu}:\mathfrak g\to \mathfrak
v(M,\mathcal O_M)$, given by the formula $X\mapsto (\varepsilon
\times \id)_2 \circ(X\oplus\, 0)\circ \mu_2$, where $X\in \mathfrak
g$. An action $\mu$ is called  {\it transitive} if the mapping
$\ev_{x} \circ \widetilde{\mu}$ is surjective for all $x\in M$, see
\cite{onipi}. In this case the supermanifold $(M,\mathcal O_M)$ is
called $(G,\mathcal O_G)$-{\it homogeneous}. A supermanifold
$(M,\mathcal O_M)$ is called {\it homogeneous}, if it possesses a
transitive action of some Lie supergroup.

Let $\phi=(\phi_1,\phi_2):(M,\mathcal O_M)\to (M_1,\mathcal
O_{M_1})$ be a morphism of supermanifolds. Denote by $(\d\!\phi)_x$
the differential of $\phi$ at $x\in M$ (see \cite{ley, man}). Let
$\dim(M,\mathcal O_M)=n|m$, $\dim(M_1,\mathcal O_{M_1})=k|l$. The
morphism $\phi:(M,\mathcal O_M)\to (M_1,\mathcal O_{M_1})$ is called
a {\it submersion at $p\in M$} if $n\geq k$, $m\geq l$ and there
exist two neighborhoods $(U,\mathcal O_M)$ of $p$ and $(V,\mathcal
O_{M_1})$ of $q=\phi_1(p)$ with the coordinates $(x_i,\,\xi_j)$ and,
respectively, $(y_s,\,\eta_t)$ such that $\phi_2|(U,\mathcal O_M)$
is given by the formulas:
\begin{equation}
\label{submersia}
\begin{array}{l}
\phi_2(y_i)=x_i,\,\, i=1,\ldots,k,\quad \phi_2(\eta_j)=\xi_j,\,\,
j=1,\ldots,l.
\end{array}
\end{equation}
 This definition is equivalent to the requirement that the mapping $(\d\!\phi)_p$ is surjective
(see \cite{ley, man}). A morphism $\phi=(\phi_1,\phi_2):(M,\mathcal
O_M)\to (M_1,\mathcal O_{M_1})$ is called an {\it  immersion at
$p\in M$}, if $n< k$, $m< l$, and there are neighborhoods
$(U,\mathcal O_M)$ and $(V,\mathcal O_{M_1})$ of $p$ and
$q=\phi_1(p)$ with coordinates $(x_i,\,\xi_j)$ and, respectively,
$(y_s,\,\eta_t)$, such that the morphism $\phi_2|(U,\mathcal O_M)$
is given by the formulas:
\begin{equation}
\label{submersia}
\begin{array}{l}
\phi_2(y_i)=x_i,\,\, i=1,\ldots,m,\quad \phi_2(y_i)=0,\,\, i>m,\\
\phi_2(\eta_j)=\xi_j,\,\, j=1,\ldots,n,\quad \phi_2(\eta_j)=0,\,\,
j>n.
\end{array}
\end{equation}
 This definition is equivalent to the requirement that the mapping $(\d\!\phi)_p$ is injective
(see \cite{ley, man}).

Let $(M,\mathcal O_{M})$ be a supermanifold. A supermanifold
$(U,\mathcal O_{M})$, where $U$ is an open subset in $M$, is called
an {\it open subsupermanifold} of $(M, \mathcal O_M)$. Suppose that
$M\subset M_1$ is a topological subspace and denote by $\phi_1:M\to
M_1$ the embedding. A supermanifold $(M,\mathcal O_{M})$ is called a
{\it subsupermanifold} of a supermanifold $(M_1,\mathcal O_{M_1})$
if there is a morphism $\phi:(M,\mathcal O_{M})\to (M_1,\mathcal
O_{M_1})$, such that the differential $(\d\!\phi)_p$ is injective at
every point $p\in M$ and the first component of $\phi$ coincides
with $\phi_1$. In this case we will use the notation $(M,\mathcal
O_{M}) \subset (M_1,\mathcal O_{M_1})$.

Let $(K,\mathcal O_K)\subset (M,\mathcal O_M)$ be a
subsupermanifold, $\mathcal{I}\subset \mathcal O_M$ the sheaf of
ideals corresponding to $(K,\mathcal O_K)$, and
$\mathcal{J}\subset\mathcal O_{G\times M}$ the sheaf of ideals
corresponding to $(G\times K,\mathcal O_{G\times K})\subset(G\times
M,\mathcal O_{G\times M})$ (see \cite{ley}). We will say that
$(K,\mathcal O_K)$ is {\it $\mu$-invariant} if the following
conditions hold:
\begin{enumerate}
  \item $\mu_1(G,K)\subset K$,
  \item $\mu_2(\mathcal{I})\subset
(\mu_1)_*\mathcal{J}$.
\end{enumerate}
A {\it Lie subsupergroup} of a Lie supergroup $(G,\mathcal O_G)$ is
a subsupermanifold $(H,\mathcal O_H)\subset (G,\mathcal O_G)$, such
that $e\in H$ and
$(H,\mathcal O_H)$ is $\nu$- and $\iota$-invariant.

If $\phi:(M,\mathcal O_M)\to (M_1,\mathcal O_{M_1})$ is a morphism
of supermanifolds and $(N,\mathcal O_N)\subset (M,\mathcal O_M)$ is
a subsupermanifold, we denote by $\phi|_{N}$ the composition:
$$
(N,\mathcal O_N) \hookrightarrow (M,\mathcal
O_M)\stackrel{\Phi}{\longrightarrow} (M_1,\mathcal O_{M_1}).
$$
Let $(G,\mathcal O_G)$ be a Lie supergroup. To each point $g\in G$
assign a morphism
$\widehat{g}=(\widehat{g}_1,\widehat{g}_2):(\pt,\Bbb C)\to
(G,\mathcal O_G)$. Namely, let $\widehat{g}_1(\pt)=g$ and define
$\widehat{g}_2: \mathcal O_G\to(\widehat{g}_1)_*(\Bbb C)$ by
$\widehat{g}_2(f_x)=0$ for $f_x\in (\mathcal O_G)_x$ if $x\ne g$ and
$\widehat{g}_2(f_g)=f_g(g)$ for $f_g\in (\mathcal O_G)_g$. Denote by
$l_g$, $g\in G$, the composition of the morphisms
$$
(G,\mathcal O_G)\xrightarrow{\sim}(\pt,\Bbb C)\times (G,\mathcal
O_G)\xrightarrow{\widehat{g}\times id} (G,\mathcal O_G)\times
(G,\mathcal O_G) \xrightarrow{\nu}(G,\mathcal O_G).
$$
and define $r_g$ in a similar way.
 Since there exists the inverse morphism $l_{g^{-1}}$,
the morphism $l_g$ is an automorphism of the supermanifold
 $(G,\mathcal O_G)$, and the same is true for $r_{g}$.
We refer the reader to \cite{ley, man} for the definition of a
superdomain and the proof of the following inverse function theorem.

\medskip

\t \label{ob obratnoj funkcii} {\sl Let $(U,\mathcal{O}_U)$ and
$(V,\mathcal{O}_V)$ be two superdomains with coordinate systems
$(x_i,\xi_j)$ and $(y_s,\eta_t)$. Let $\phi: (U,\mathcal{O}_U)\to
(V,\mathcal{O}_V)$ be a morphism and let $u\in U$. Then the
following conditions are equivalent:
\begin{itemize}
  \item $\phi$ is an isomorphism in some neighborhood of $u$;
  \item $(\d\!\phi)_u$ is an isomorphism.
\end{itemize}
}\noindent Given supermanifolds $(M_i,\mathcal{O}_{M_i})$,
$i=1,\ldots,n$, we will denote by $\pr_{M_i}^{M_1\times \cdots\times
M_n}$ the natural projection $(M_1,\mathcal{O}_{M_1})\times
\cdots\times (M_n,\mathcal{O}_{M_n})\to (M_i,\mathcal{O}_{M_i})$.

\bigskip

\bigskip

\begin{center}
{\bf 2. The structure of a supermanifold on $G/H$}
\end{center}

\bigskip

Let $(K,\mathcal O_{K})\subset (U,\mathcal O_{G}|_U)$ be a
subsupermanifold of an open subsupermanifold $(U,\mathcal O_{G}|_U)$
in a Lie supergroup $(G,\mathcal O_{G})$, $\mathcal{I}_K$ the
corresponding sheaf of ideals and $\varphi$ is an isomorphism of
$(G,\mathcal O_{G})$. We will denote by $(\varphi (K),\mathcal
O_{\varphi (K)})$ the subsupermanifold of
 $(\varphi (U),\mathcal
O_{G}|_{gU})$, where $\varphi (K):=\varphi_1(K)$, $ \mathcal
O_{\varphi (K)}:=(\mathcal O_{G}|_{\varphi (U)})/
(\varphi^{-1})_2(\mathcal{I}_K). $ Sometimes we will use the
notation $(gK,\mathcal O_{gK})$ for $(l_g(K),\mathcal O_{l_g(K)})$.
The subsupermanifold $(Kg,\mathcal O_{Kg})$ is defined analogously.
The following proposition is well known.

\medskip

\l\label{obratn_morphizm} {\sl Let $\varphi=(\varphi_1,\varphi_2):
(M,\mathcal O_M)\to (N,\mathcal O_N)$ be a  morphism of
supermanifolds. Assume that $\varphi_1:M\to N$ is a homeomorphism
and $\varphi$ is a local isomorphism. Then $\varphi$ is an
isomorphism. $\Box$}

\medskip

Let $(N,\mathcal O_{N})$ and $(S,\mathcal O_{S})$ be two
subsupermanifolds of $(M,\mathcal O_{M})$. The subsupermanifold
$(S,\mathcal O_{S})$ is called {\it transversal to the
subsupermanifold $(N,\mathcal O_{N})$ at a point $x\in S\cap N$}, if
$T_x(M,\mathcal O_{M})=T_x(N,\mathcal O_{N})\oplus T_x(S,\mathcal
O_{S})$.

\medskip

\t\label{triv rassl} {\sl Let $(G,\mathcal O_G)$ be a Lie supergroup
and $(H,\mathcal O_H)$ be a Lie subsupergroup of $(G,\mathcal O_G)$.
Suppose that $(S',\mathcal O_{S'})$ is a transversal
subsupermanifold to $(H,\mathcal O_H)$ at the point $e$ (the
identity element of $G$). Then there is a subsupermanifold
$(S,\mathcal O_S)$, such that
 $\nu|_{S\times H} :(S,\mathcal
O_S)\times (H,\mathcal O_H)\to (G,\mathcal O_G)$ is an isomorphism
of
 $(S,\mathcal O_S)\times (H,\mathcal
O_H)$  onto an open subsupermanifold $(U,\mathcal
O_G)\subset(G,\mathcal O_G)$. }

\medskip

\noindent{\it Proof.} Let $\dim(G,\mathcal O_G)=n|m$. First we will
show that there exists a transversal subsupermanifold to
$(H,\mathcal O_H)$ at $e$. From the definition of a subsupermanifold
it follows that there is a superdomain $(V,\mathcal O_V)\subset
(G,\mathcal O_G)$ containing $e$, such that the subsupermanifold
$(H,\mathcal O_H)$ is given by the equations $x_i=0$, $i\in
\Gamma_1$, and $\xi_j=0$, $j\in \Gamma_2$, where $\Gamma_1\subset
\{1,\ldots,n\}$, $\Gamma_2\subset \{1,\ldots,m\}$. Denote by
$(S',\mathcal O_{S'})$ the subsupermanifold in $(V,\mathcal O_V)$
given by the equations $x_i=0$, $i\notin \Gamma_1$, and $\xi_j=0$,
$j\notin \Gamma_2$. Obviously, $(S',\mathcal O_{S'})$ is transversal
to $(H,\mathcal O_H)$.

Since $\nu|_{S\times H}\circ (\varepsilon\times \id)=\id$ and
$\nu|_{S\times H}\circ (\id\times \varepsilon)=\id$, it follows that
the differential $(\d\nu|_{S\times H})_{(e,e)}:T_{e} (S_1,\mathcal
O_{S_1})\oplus T_{e}(H,\mathcal O_H)=T_{e}(G,\mathcal O_G) \to
T_{e}(G,\mathcal O_G)$ is precisely the mapping
\begin{equation}
\label{diff}
 (v_s, v_h)\mapsto v_s+v_h.
\end{equation}
 Now it is easy to see that
$(\d\nu|_{S\times H})_{(e,e)}$ is an isomorphism. From Theorem
\ref{ob obratnoj funkcii}, it follows that there are superdomains
$(S,\mathcal O_{S})\subset (S',\mathcal O_{S'})$, $(H',\mathcal
O_{H})\subset (H,\mathcal O_H)$ and $(W,\mathcal O_G)\subset
(G,\mathcal O_G)$, such that $\nu|_{S\times H'}:(S,\mathcal
O_{S})\times (H',\mathcal O_{H})\to (W,\mathcal O_G)$ is an
isomorphism.

Fix a point $h \in H$. By the associativity axiom for Lie
supergroups we have the following commutative diagram:
$$
\begin{CD}
(S,\mathcal O_{S})\times (H',\mathcal O_{H})@>{id\times r_h|_{H}}>>
(S,\mathcal O_{S})\times (H'h,\mathcal O_{H})
\\
@V{\nu|_{S\times H'}}VV @VV{\nu|_{S\times H'h}}V\\
(G,\mathcal O_G)@>{r_h}>>  (G,\mathcal O_G)
\end{CD}.
$$
In other words, the morphism $\nu $, restricted to $(S,\mathcal
O_{S})\times (H'h,\mathcal O_{H})$,
 is equal to $ r_h\circ
\nu|_{S\times H'}\circ (id\times r_h|_H)^{-1}$. We have
$$
(S,\mathcal O_{S})\times (H,\mathcal O_{H})=\bigcup_{h\in H}\bigl
\{(S,\mathcal O_{S})\times (H'h,\mathcal O_{H})\bigr \}.
$$
Therefore $\nu|_{S\times H}$ is a local isomorphism. By a well known
argument from the geometric theory of homogeneous spaces we may
assume that $(\nu|_{S\times H})_1:S\times H\to G$ is a
homeomorphism. Using Lemma \ref{obratn_morphizm}, we get that
 $\nu|_{S\times H}:(S,\mathcal O_S)\times (H,\mathcal O_H)\to (U,\mathcal O_U)$ is an isomorphism.
$\Box$
\medskip

Now we will give the definition of a supermanifold with the
underlying
 manifold $G/ H$, corresponding to a Lie
supergroup $(G,\mathcal O_G)$ and a subsupergroup $(H,\mathcal
O_H)\subset(G,\mathcal O_G)$. Let $p_1:G\to G/H$, $g\mapsto gH$ be
the natural mapping. Fix a transversal subsupermanifold $(S,\mathcal
O_{S})\subset (G,\mathcal O_{G})$
 to $(H,\mathcal O_{H})$ at the point $e$ so that $\nu|_{S\times H}$
 is an isomorphism (see Theorem \ref{triv rassl}).
 The mapping $p_1$ maps
$S$ homeomorphically  onto a domain $V\subset G/H$. Denote by
$\mathcal{O}_V$ the sheaf $(p_1)_*(\mathcal O_S)$ and identify
$(S,\mathcal O_{S})$ with $(V,\mathcal O_{V})$. Recall that the open
subsupermanifold $(U,\mathcal O_G)\subset (G,\mathcal O_G)$ is
defined in Theorem \ref{triv rassl}. Let $p_U=(p_1,(p_U)_2)$ be such
a morphism that the following diagram is commutative:
$$
\begin{CD}
(S,\mathcal O_{S})\times (H,\mathcal O_{H})@>{\nu|_{S\times H}}>>
(U,\mathcal O_{G})
\\
@V{\pr_V^{S\times H}}VV @VV{p_U}V\\
(V,\mathcal O_V)@= (V,\mathcal O_{V})
\end{CD}\,\,\,\,\,\,\,\,\,.
$$
 We can define $(gV,\mathcal O_{gV})$ analogously. Namely, it is
clear that $\nu|_{gS\times H}:(gS,\mathcal O_{gS})\times (H,\mathcal
O_H)\to (gU,\mathcal O_{gU})$ is an isomorphism and $p_1:gS\to gV$
is a homeomorphism. Denote by $\mathcal O_{gV}$ the sheaf
$(p_1)_*(\mathcal O_{gS})$ on $gV$ and identify the subsupermanifold
$(gS,\mathcal O_{gS})$ with $(gV,\mathcal O_{gV})$. Define a
morphism $p_{gU}:
 (gU,\mathcal O_{G})\to (gV,\mathcal O_{gV})$ by
$$p_{gU}= \pr_{gV}^{gS\times H}\circ (\nu|_{gS\times H})^{-1}.
$$
Now we need the following lemma.

\medskip
\l\label{funkc na triv rassloen} {\it Let $(W,\mathcal O_W)$ be a
superdomain, $(F,\mathcal O_F)$ a supermanifold, and $(M,\mathcal
O_M)= (W,\mathcal O_W)\times (F,\mathcal O_F)$. Denote by $(x_i)$
all, i.e. even and odd, coordinates on $(W,\mathcal O_W)$. Then any
function $f\in \mathcal O_M(M)$ can be written in the form $f=\sum_j
h_j g_j$, where $(h_j)$ is a maximal independent system of
polynomials in $x_i$ and $g_j\in\mathcal O_F(F)$.}

\medskip

\noindent{\it Proof.} Fix a coordinate neighborhood $(N,\mathcal
O_F)$ of the supermanifold $(F,\mathcal O_F)$. The function
$f|_{W\times N}$ has the form $\sum_jh_jt_j^N$, where $t_j^N\in
\mathcal O_F(N)$. Let $(\tilde{N},\mathcal O_F)$ be another
coordinate neighborhood of $(F,\mathcal O_F)$, such that $N\cap
\tilde N\ne \varnothing$. Then the function $f|_{W\times \tilde N}$
has the form
$$
f=\sum_j  h_j t_{j}^{\tilde N},
$$
where $t_{j}^{\tilde N}\in \mathcal O_F(\tilde N)$. It is obvious
that $t_j^N=\tilde t_{j}^{\tilde N}$ in $N\cap\tilde{N}$. Now we
choose an atlas $\{ (N,\mathcal O_F)\}$ on the supermanifold
$(F,\mathcal O_F)$. The functions $t_{j}^N$ are holomorphic in all
coordinate neighborhoods of the chosen atlas and coincide on their
intersections. It follows that there are $g_{j}\in H^0(F,\mathcal
O_F)$ such that $g_j|_{N}=t_j^N$, and so we can write $f=\sum_j h_j
g_j$.$\Box$

\medskip

Let $W\subset G/H$ be an open set. A function $f\in
\mathcal{O}_G(p_1^{-1}(W))$ is called {\it $(H,\mathcal
O_{H})$-right invariant} if $(\nu|_{G\times
H})_2(f)=(\pr^{p_1^{-1}(W)\times H}_{p_1^{-1}(W)})_2(f)$ (see
\cite{Kostant}).

\medskip
\l\label{f-(H,O_H)-invariantna} {\it Let $W$ be an open set in $gV$,
$f\in \mathcal{O}_G(p_{gU}^{-1}(W))$. Then $f\in \Im(p_{gU})_2$ if
and only if $f$ is a $(H,\mathcal O_{H})$-right invariant function.}

\medskip
\noindent{\it Proof.}  By construction, we have $\Im
(p_{gU})_2=((\nu|_{gS\times H})^{-1})_2((\pr_{gS}^{gS\times
H})_2(\mathcal O_{gS}))$. By associativity of the multiplication in
$(G,\mathcal O_{G})$, the following diagram is commutative:
\begin{equation}
\label{assoc}
\begin{CD}
(gU, \mathcal O_{G})\times (H,\mathcal
O_{H})@>{\nu|_{G\times H}}>>(gU, \mathcal O_{G})\\
@A{\nu|_{gS\times H}\times \id}AA @AA{
\nu|_{gS\times H}.}A\\
(gS,\mathcal O_{gS})\times (H,\mathcal O_{H})\times (H,\mathcal
O_{H})@>{\id\times \nu|_{H\times H}}>>(gS,\mathcal O_{gS})\times (H,\mathcal O_{H})\\
\end{CD}
\end{equation}
Now we can see that if $f\in \Im(p_{gU})_2$ then $f$ is a
$(H,\mathcal O_{H})$-right invariant function. Indeed, in this case
$(\nu|_{gS\times H})_2(f)\in (\pr_{gS}^{gS\times H})_2(\mathcal
O_{gS})$ and
$$
(\id_2\times (\nu|_{H\times H})_2)((\nu|_{gS\times
H})_2(f))=(\pr_{gS}^{gS\times H\times H})_2((\nu|_{gS\times
H})_2(f)).
$$
Applying the isomorphism  $(\nu|_{gS\times H}\times \id)^{-1}$ to
the right hand side, we get:
$$
(\nu|_{gS\times H}\times \id)^{-1}_2\circ(\pr_{gS}^{gS\times H\times
H})_2((\nu|_{gS\times H})_2(f)) =(\pr_{gU}^{gU \times H})_2(f).
$$
By the commutativity of (\ref{assoc}), we obtain:
$$
(\nu|_{G\times H})_2(f)=(\pr_{p_1^{-1}(W)}^{p_1^{-1}(W) \times
H})_2(f).
$$

Conversely, let $f\in \mathcal{O}_G(p_{gU}^{-1}(W))$ be a
$(H,\mathcal O_{H})$-right invariant function. Without loss of
generality assume that $(gS,\mathcal O_{gS})$ is a coordinate
neighborhood. In the present setting, there will be no confusion to
denote the even and odd coordinates by the same letters $(x_i)$. By
Lemma \ref{funkc na triv rassloen}, we have
$$
(\nu|_{gS\times H})_2(f)=\sum h_ig_i,
$$
where $(h_i)$ is a maximal independent system of monomials in $x_i$,
$g_i\in \mathcal O_{H}(H)$. By the definition of a $(H,\mathcal
O_{H})$-right invariant function and the commutativity of
(\ref{assoc}), we get
$$
(\id\times \nu|_{H\times H})_2(\sum h_ig_i)=(\pr_{gS\times
H}^{gS\times H\times H})_2(\sum h_ig_i)= \sum (\pr_{gS}^{gS\times
H\times H})_2(h_i)  (\pr_{H_1}^{gS\times H\times H})_2(g_i),
$$
where $\pr_{H_1}^{gS\times H\times H}$ is the projection onto the
second factor. On the other hand,
$$ (\id\times \nu|_{H\times H})_2(\sum h_ig_i)= \sum
(\pr_{gS}^{gS\times H\times H})_2(h_i) (\pr_{H\times H}^{gS\times
H\times H})_2 ((\nu|_{H\times H})_2(g_i)).
$$
The independence of $(h_i)$ implies that $(\pr_{H_1}^{ H\times
H})_2(g_i)= (\nu|_{H\times H})_2(g_i)$, where $\pr_{H_1}^{H\times
H}$ is the projection onto the first factor. Equivalently,
 the functions $g_i$ are $(H,\mathcal O_{H})$-right
invariant. We have reduced our assertion to the following one.

\medskip

\noindent $(*)$ {\it If a function $g\in \mathcal O_{H}(H)$ is
$(H,\mathcal O_{H})$-right invariant, then $g=const$.}

\smallskip
\noindent{\it Proof of $(*)$.} As above, let $\pr^{H\times
H}_{H_1}:(H,\mathcal O_{H})\times (H,\mathcal O_{H})\to (H,\mathcal
O_{H})$ denote the projection onto the first factor. By the
definition of a $(H,\mathcal O_{H})$-right invariant function, we
have $(\pr^{H\times H}_{H_1})_2(g)=(\nu|_{H\times H})_2(g)$. Now we
get
$$
g(e)=(\varepsilon \times \id)_2((\pr^{H\times
H}_{H_1})_2(g))=(\varepsilon \times \id)_2((\nu|_{H\times
H})_2(g))=g,
$$
where the last equality follows from the identity axiom $\nu\circ
(\varepsilon \times \id)=\id$. Therefore $g= g(e)$ showing $(*)$.
This completes the proof of Lemma \ref{f-(H,O_H)-invariantna}.
$\Box$

\medskip

\t\label{U to gU} {\sl The charts $(gV,\mathcal O_{gV})$ constitute
a holomorphic atlas on $G/H$. }

\medskip

\noindent{\it Proof.} Suppose that $g_1V\cap g_2V\ne \emptyset$. Let
us prove that there is a morphism $\Psi_{g_1V,g_2V}:(g_1V\cap
g_2V,\mathcal O_{g_1V}) \to (g_1V\cap g_2V,\mathcal O_{g_2V})$ such
that $p_{g_2U}=\Psi_{g_1V,g_2V} \circ p_{g_1U}$. Obviously,
$(\Psi_{g_1V,g_2V})_1=\id$. Let us define the second component
 $(\Psi_{g_1V,g_2V})_2$.

If $f\in \mathcal O_{g_2V}|_{g_1V\cap g_2V}$ then
$(p_{g_2U})_2(f)\in \mathcal O_{G}|_{g_1U\cap g_2U}$. By Lemma
\ref{f-(H,O_H)-invariantna}, the function $(p_{g_2U})_2(f)$ is
$(H,\mathcal O_{H})$-right invariant. Therefore $(p_{g_2U})_2(f)\in
\Im ((p_{g_1U})_2)$. By construction, the map $(p_{g_1U})_2$ is
injective, and so there is a unique function $g\in \mathcal
O_{g_1V}|_{g_1V\cap g_2V}$, such that
$(p_{g_1U})_2(g)=(p_{g_2U})_2(f)$. We put
$(\Psi_{g_1V,g_2V})_2(f):=g$. The cocycle condition is obviously
fulfilled.$\Box$

\medskip
We denote by $(G/H,\mathcal O_{G/H})$ the supermanifold defined by
the holomorphic atlas constructed above. From the definition of
transition functions between the charts $(gV,\mathcal O_{gV})$, we
get a morphism $p:(G,\mathcal O_{G})\to (G/H,\mathcal O_{G/H})$ with
$p|_{gU}=p_{gU}$ for all $g\in G$. Now we can generalize Lemma
\ref{f-(H,O_H)-invariantna}.

\medskip
\l\label{f-(H,O_H)-invariantna2} {\it Let $W\subset G/H$ be an open
set and  $f\in \mathcal{O}_G(p_1^{-1}(W))$. Then $f\in \Im(p_2)$ if
and only if $f$ is $(H,\mathcal O_{H})$-right invariant.} $\Box$

\medskip


We now define an action of the Lie supergroup $(G,\mathcal O_{G})$
on the supermanifold $(G/H,\mathcal O_{G/H})$, which will be used in
the proof of the main result of this section (Theorem
\ref{struktura_supera}). Denote by $i_{g}:(gV,\mathcal
O_{G/H})=(gS,\mathcal O_{gS})\to (gU,\mathcal O_{G})$ the natural
embedding. 
 Let
$\alpha_{gV}$ be the following composition:
\begin{equation}\label{deydef}
\begin{CD}
(G,\mathcal O_{G})\times (gV,\mathcal O_{G/H})@>{\id\times
i_g}>>(G,\mathcal O_{G})\times (gU,\mathcal
O_{G})@>{\nu}>>(G,\mathcal O_{G}) @>{p}>>(G/H,\mathcal O_{G/H})
\end{CD}.
\end{equation}
Suppose $g_1V\cap g_2V\ne \emptyset$. We claim that
$\alpha_{g_1V}|_{G\times g_1V\cap g_2V}=\alpha_{g_2V}|_{G\times
g_1V\cap g_2V}$. Obviously,\linebreak $(\alpha_{g_1V})_1|_{G\times
(g_1V\cap g_2V)}=(\alpha_{g_2V})_1|_{G\times (g_1V\cap g_2V)}$, and
so we only have to prove our equality for the second compo\-nents of
the morphisms $\alpha_{g_1V}$ and $\alpha_{g_2V}$. Let $W\subset
G/H$ be an open set and $f\in \mathcal O_{G/H}(W)$. It suffices to
check that
\begin{equation}\label{alpha_loc}
(\id_2\times (i_{g_1})_2)(\nu_2\circ p_2(f))=(\id_2\times
(i_{g_2})_2)(\nu_2\circ p_2(f)).
\end{equation}
From the
associativity of the multiplication in $(G,\mathcal{O}_G)$ and from
Lemma \ref{f-(H,O_H)-invariantna2}, it follows that
$$
(\id_2\times (\nu|_{G\times H})_2)(\nu_2\circ p_2(f))=(\id_2\times
(\pr^{G\times H}_G)_2)(\nu_2\circ p_2(f)).
$$
Let $(O,\mathcal O_{G})$ and $(K,\mathcal O_{G})$ be such open sets
in $(G,\mathcal O_{G})$ that $O\times K\subset
\nu_1^{-1}(p_1^{-1}(W))$. We assume that $(O,\mathcal O_{G})$ is a
coordinate neighborhood with even and odd coordinates $(y_j)$ and
$K=p_1^{-1}(N)$, where $N$ is an open set in $G/H$. Denote by
$(t_i)$ a maximal independent system of monomials in $y_j$. By Lemma
\ref{funkc na triv rassloen} we can write the function $\nu_2\circ
p_2(f)|_{O\times K}$ in the form $\sum t_i s_i$, where $s_i\in
\mathcal O_{G}(K)$. We have in $O\times K\times H$
$$
\begin{array}{l}
\sum t_i (\pr^{G\times H}_G)_2(s_i)=(\id_2\times (\pr^{G\times
H}_G)_2)(\sum t_i s_i)= (\id_2\times (\pr^{G\times
H}_G)_2)(\nu_2\circ p_2(f)) = \\
\\
 =(\id_2\times (\nu|_{G\times H})_2)(\nu_2\circ
p_2(f))=(\id_2\times (\nu|_{G\times H})_2)(\sum t_i s_i)=\sum t_i
(\nu|_{G\times H})_2(s_i).
\end{array}
$$
Now the equality $(\pr^{G\times H}_G)_2(s_i)= (\nu|_{G\times
H})_2(s_i)$ follows from the independence of the system $(t_i)$. In
other words, we get that the functions $s_i$ are $(H,\mathcal
O_{H})$-right invariant. By Lemma \ref{f-(H,O_H)-invariantna2} there
are functions $h_i\in \mathcal O_{G/H}(N)$, such that
$p_2(h_i)=s_i$. Moreover,
$$
(\id_2\times (i_{g_1})_2)(\nu_2\circ p_2(f)|_{O\times
K})=(\id_2\times (i_{g_1})_2)(\sum t_i s_i)=\!\sum t_i
(i_{g_1})_2(s_i)=\! \sum t_i (i_{g_1})_2(p_2(h_i)).
$$
Now we use the trivial equality $p\circ i_{g}=\id$ for all $g\in G$.
$$
\sum t_i (i_{g_1})_2(p_2(h_i))=\sum t_i h_i.
$$
Similarly, we obtain
$$
(\id_2\times (i_{g_2})_2)(\nu_2\circ p_2(f)|_{O\times K})= \sum t_i
h_i.
$$
So we have shown that
$$
(\id_2\times (i_{g_1})_2)(\nu_2\circ p_2(f))|_{O\times
N}=(\id_2\times (i_{g_2})_2)(\nu_2\circ p_2(f)) |_{O\times N}.
$$
This implies (\ref{alpha_loc}), and so we get a morphism $\alpha$,
such that $\alpha|_{G\times gV}=\alpha_{gV}$. Clearly, $\alpha$ is
an action on the supermanifold $G/H$. We have proved the following
theorem.

\medskip

\t\label{struktura_supera} {\sl  There exists a supermanifold
$(G/H,\mathcal O_{G/H})$, such that the natural action of $G$ on
$G/H$ induces a transitive action of $(G,\mathcal{O}_G)$ on
$(G/H,\mathcal O_{G/H})$. The action of the Lie supergroup
$(G,\mathcal O_{G})$ on $(G/H,\mathcal O_{G/H})$ is given by
(\ref{deydef}).}$\Box$
\medskip

\bigskip

\bigskip

\begin{center}
{\bf 3. Stationary Lie subsupergroup}
\end{center}

\bigskip

Let $\mu:(G,\mathcal O_G)\times (M,\mathcal O_M)\to (M,\mathcal
O_M)$ be an action of a Lie supergroup $(G,\mathcal O_G)$ on a
supermanifold $(M,\mathcal O_M)$ and let $\mathfrak{g}$ be the Lie
superalgebra of the Lie supergroup $(G,\mathcal O_G)$. Denote by
$\mu_x:(G,\mathcal O_G)\to (M,\mathcal O_M)$, $x\in M$, the
composition of morphisms
$$
(G,\mathcal O_G)\times (\pt,\mathbb C)\xrightarrow{id\times
\widehat{x}}(G,\mathcal O_G)\times (M,\mathcal
O_M)\xrightarrow{\mu}(M,\mathcal O_M),
$$
where $\widehat{x}=(\widehat{x}_1,\widehat{x}_2):(\pt,\mathbb C)\to
(M,\mathcal O_M)$, $\widehat{x}_1(\pt)=x$, $\widehat{x}_2(f)=f(x)$,
$f\in \mathcal O_M$. Also, let $\overline{l}_g:(M,\mathcal O_M)\to
(M,\mathcal O_M)$, $g\in G$, be the composition of morphisms
$$
(M,\mathcal O_{M})\simeq(\pt,\mathbb C)\times (M,\mathcal
O_{M})\xrightarrow{\widehat{g}\times id}(G,\mathcal O_G)\times
(M,\mathcal O_M)\xrightarrow{\mu}(M,\mathcal O_M),
$$
where $\widehat{g}$ was defined in Section 1.

\medskip

\l\label{dlya glupyh} {\it  We have $\ev_x\circ
\widetilde{\mu}(X)=(\d\!\mu_{x})_{e}(\ev_e(X))$, $X\in
\mathfrak{g}$.
 The action $\mu$ is
transitive if and only if $\mu_x$ is a submersion at $e\in G$ for
all $x\in M$. }

\medskip

\noindent{\it Proof.} The second assertion follows from the first
one. Let us prove the first assertion. By definition we get
$$
\ev_x(\widetilde{\mu}(X))(f)=(\widetilde{\mu}(X)(f))(x),\,\,\,\,
(\d\!\mu_{x})_{e}(\ev_e(X))(f)=(X((\mu_x)_2(f)))(e)
$$
for all $f\in (\mathcal O_M)_x$. Therefore,
$$
\begin{array}{rll}
\ev_x(\widetilde{\mu}(X))(f)=&(\id_2\times
\widehat{x}_2)((\varepsilon_2\times \id_2)((X\oplus 0)\circ
\mu_2(f)))=&\\ &(\varepsilon_2\times \id_2)((\id_2\times
\widehat{x}_2)((X\oplus 0)\circ \mu_2(f)))=&\\
&\varepsilon_2(X((\id_2\times \widehat{x}_2)\circ\mu_2(f)))=
\varepsilon_2(X((\mu_x)_2(f)))=&(\d\!\mu_{x})_{e}(X_{e})(f).\Box
\end{array}
$$

\smallskip

\noindent The following lemma is a consequence of the axioms of
action.

\medskip

\l\label{sdvig_deystvija} {\sl We have $\mu_x\circ r_g=\mu_{gx}$,
$\mu_x\circ l_g=\overline{l}_g\circ\mu_{x}$. $\Box$}
\medskip

 Consider a superdomain
$(V,\mathcal O_M)\subset (M,\mathcal O_M)$, such that $x$ is
contained in $V$, with even coordinates $(y_i)$ and odd coordinates
 $(\eta_j)$. Suppose that $x$ has the coordinates $y_i=0$,
$\eta_j=0$. Denote by $(U,\mathcal O_G)$ a superdomain in
$((\mu_x)_1^{-1}(V),\mathcal O_G)$ with even coordinates $(x_s)$ and
odd coordinates $(\xi_t)$. Let the morphism $\mu_x|(U,\mathcal O_G)$
be given by the equations:
$$
(\mu_x)_2(y_i)=\phi^i_U(x_s,\xi_t), \quad
(\mu_x)_2(\eta_j)=\psi^j_U(x_s,\xi_t).
$$
Denote by $\mathcal{I}_U$ the sheaf of ideals in the structure sheaf
of $U$ defined as follows: if $U\cap (\mu_x)_1^{-1}(V)\ne\emptyset$
then $\mathcal{I}_U$ is generated by the functions
$\phi^i_U(x_s,\xi_t)$, $\psi^j_U(x_s,\xi_t)$, otherwise
$\mathcal{I}_U:=\mathcal O_G|_U$. The sheafs $\mathcal{I}_{U_1}$ and
$\mathcal{I}_{U_2}$ coincide on the intersection $U_1\cap U_2$.
Therefore there is a sheaf of ideals $\mathcal{I}$, such that
$\mathcal{I}|_U=\mathcal{I}_U$. As usual, denote by $G_x$ the
stationary subgroup of the action $\mu _1$ at $x$ and consider the
ringed space $(G_x,\mathcal O_{G_x})$, where  $\mathcal
O_{G_x}=(\mathcal O_G/\mathcal{I})|_{G_x}$.

If $\mu$ is transitive then $(G_x,\mathcal O_{G_x})$ is a
subsupermanifold of $(G,\mathcal O_{G})$. Indeed, by Lemma \ref{dlya
glupyh} the morphism $\mu_x$ is a submersion at $e$. From Lemma
\ref{sdvig_deystvija} it follows that $\mu_x$ is a submersion at
every point $g\in G$. By the definition of a submersion, there is a
coordinate neighborhood $(W,\mathcal O_{G})$ of $g\in G_x$ with
coordinates $(x_i;\xi_j)$, $i=1,\ldots,p$, $j=1,\ldots,q$ and a
coordinate neighborhood of $x$ with coordinates $(y_a;\eta_b)$,
$a=1,\ldots,k$, $b=1,\ldots,l$,  such that $\mu_x|(W,\mathcal
O_{G})$ is given by the following formulas:
\begin{equation}
\label{mu_m}
\begin{split}
(\mu_x)_2(y_a)=x_a,\,\,\,a=1,\ldots,k,\quad
(\mu_x)_2(\eta_b)=\xi_b,\,\,\,b=1,\ldots,l.
\end{split}
\end{equation}
Without loss of generality assume that the point $x$ is given by the
system of equations $y_a=0$, $a=1,\ldots,k$, $\eta_b=0$,
$b=1,\ldots,l$.
 Then $(G_x\cap W,\mathcal O_{G_x})$ is isomorphic to
a superdomain with coordinates $x_i, \,i=k+1,\ldots,p$,
$\xi_j,\,j=l+1,\ldots,q$.


It is not hard to prove that $\iota^2=\id$. Indeed,
$$
\nu\circ (\id\times \nu)\circ (\iota^2,\iota,\id)= \nu \circ
(\iota^2,\nu \circ (\iota,\id))=\nu \circ (\iota^2,\varepsilon)
=\iota^2
$$
and
$$
\nu\circ (\id\times \nu)\circ (\iota^2,\iota,\id)=\nu \circ
(\nu\times \id)\circ ((\iota,\id)\times \id)\circ
(\iota,\id)=\nu\circ (\varepsilon\times \id)\circ (\iota,\id) =\id.
$$
We will use this conclusion in the proof of the following theorem.

\medskip
\t\label{G_mgrupp operacii} {\it Suppose $(G_x,\mathcal
O_{G_x})\subset (G,\mathcal O_{G})$ is a subsupermanifold. Then
$(G_x,\mathcal O_{G_x})$ is a Lie subsupergroup  of $(G,\mathcal
O_{G})$.}

\medskip

\noindent{\it Proof.} We must show that $(G_x,\mathcal O_{G_x})$ is
$\nu$-invariant and $\iota$-invariant. We check first that
 \begin{equation}
 \label{nu-invariant}
 (G_x,\mathcal O_{G_x})\ \ \text{\it is
 $\nu$-invariant}.
 \end{equation}

 \noindent Obviously, $\nu_1(G_x,G_x)=G_x$. Denote by $\mathcal{J}$
  the sheaf of ideals corresponding to the subsupermanifold $(G_x,\mathcal
O_{G_x})\times (G_x,\mathcal O_{G_x})$ of the supermanifold
$(G,\mathcal O_{G})\times (G,\mathcal O_{G})$. We have to prove that
$\nu_2(\mathcal{I})\subset \mathcal{J}$. The functions
$\phi^i_U(x_s,\xi_t)$, $\psi^j_U(x_s,\xi_t)$ generate the ideal
sheaf $\mathcal{I}|_U$ and $\mathcal{I}_y\ne (\mathcal O_{G})_y$
only for $y\in G_x\subset W:=(\mu_x)_1^{-1}(V)$. Therefore it is
sufficient to prove that $\nu_2(\phi^i_U(x_s,\xi_t))|_{W\times
W}\subset \mathcal{J}|_{W\times W}$ and
$\nu_2(\psi^j_U(x_s,\xi_t))|_{W\times W}\subset
\mathcal{J}|_{W\times W}$. By the definition of the functions
$\phi^i_U(x_s,\xi_t)$ and $\psi^j_U(x_s,\xi_t)$, we get
$\nu_2(\phi^i_U(x_s,\xi_t))=\nu_2( (\mu_x)_2(y_i))$,
$\nu_2(\psi^j_U(x_s,\xi_t))=\nu_2( (\mu_x)_2(\eta_j))$. From the
 axioms of action it follows that $\mu_x\circ \nu=\mu\circ (\id\times
\mu_x)$. Thus it suffices to prove that $(\id_2\times
(\mu_x)_2)\circ \mu_2(y_i)|_{W\times W}\subset \mathcal{J}|_{W\times
W}$ and $(\id_2\times (\mu_x)_2)\circ \mu_2(\eta_j)|_{W\times
W}\subset \mathcal{J}|_{W\times W}$.

Denote by $\pr_i:(G,\mathcal O_{G})\times (G,\mathcal O_{G})\to
(G,\mathcal O_{G})$ the projection on the $i$-th factor and by
$\widetilde{(\pr_i)_2(\mathcal{I})}$ the sheaf of ideals generated
by $(\pr_i)_2(\mathcal{I})$, $i=1,2$. We have
$\mathcal{J}=\widetilde{(\pr_1)_2(\mathcal{I})}+\widetilde{(\pr_2)_2(\mathcal{I})}$.
Further, $[(\id_2\times (\mu_x)_2)\circ \mu_2(y_i)]|_{W\times W}=
(\id_2\times (\mu_x)_2)[ \mu_2(y_i)|_{W\times V}]$. Using the
definition of $\mathcal{I}$, we get:
$$
(\id_2\times (\mu_x)_2)[ \mu_2(y_i)|_{W\times
V}]+\widetilde{(\pr_2)_2(\mathcal{I})}|_{W\times
W}=(\pr_1)_2((\mu_x)_2(y_i))|_{W\times W}+
\widetilde{(\pr_2)_2(\mathcal{I})}|_{W\times W}.
$$
Now, from $(\pr_1)_2((\mu_x)_2(y_i))|_{W\times W}\in
\widetilde{(\pr_1)_2(\mathcal{I})}|_{W\times W}$ it follows that
$$
(\id_2\times (\mu_x)_2)\circ \mu_2(y_i)|_{W\times W}\in
\widetilde{(\pr_1)_2(\mathcal{I})}|_{W\times
W}+\widetilde{(\pr_2)_2(\mathcal{I})}|_{W\times
W}=\mathcal{J}|_{W\times W}
$$
and, by the same argument, $(\id_2\times (\mu_x)_2)\circ
\mu_2(\eta_j)|_{W\times W}\subset \mathcal{J}|_{W\times W}$. This
completes the proof of (\ref{nu-invariant}).

It remains to check that
$$ 
 (G_x,\mathcal O_{G_x})\ \ \text{\it is
$\iota$-invariant}.
$$ 
Since the inclusion $\iota_1(G_x)\subset G_x$ is obvious,  we must
prove that $\iota_2(\mathcal{I})\subset\mathcal{I}$ or, in terms of
generators, $\iota_2(\phi^i_U(x_s,\xi_t))\in \mathcal{I}$ and
$\iota_2(\psi^i_U(x_s,\xi_t))\in \mathcal{I}$.

By definition of the supermanifold $(G_x,\mathcal O_{G_x})$, the
following diagram is commutative:
$$
\begin{CD}
(G_x,\mathcal O_{G_x})@>{\pr_x}>> (x,\Bbb C)\\
@| @VV{\widehat{x}}V\\
(G_x,\mathcal O_{G_x})@>{\mu_x}>>(M,\mathcal O_{M})
\end{CD}.
$$
We will rather use the commutativity of the next diagram:
\begin{equation}
\label{iota(G_x)}
\begin{CD}
(\iota_1(G_x),\mathcal O_{\iota_1(G_x)})@>{\pr_x}>> (x,\Bbb C)\\
@| @VV{\widehat{x}}V\\
(\iota_1(G_x),\mathcal O_{\iota_1(G_x)})@>{\mu_x}>>(M,\mathcal
O_{M})
\end{CD}.
\end{equation}
To show that (\ref{iota(G_x)}) is commutative, note that, by the
inverse element axiom of a Lie supergroup and by the axioms of an
action, $\mu\circ (\iota,\mu_x)= \widehat{x}\circ \pr_x$ and, in
particular, $\mu\circ (\iota,\mu_x)|_{G_x}= (\widehat{x}\circ
\pr_x)|_{G_x}$. Using the equality $\mu_x|_{G_x}=\widehat{x}$ and
the definition of the morphism $\mu_x$, we obtain the commutative
diagram:
$$
\begin{CD}
(G_x,\mathcal O_{G_x})@>{\pr_x}>> (x,\Bbb C)\\
@V{(\iota,\widehat{x})}VV @VV{\widehat{x}}V\\
(\iota_1(G_x),\mathcal O_{\iota_1(G_x)})\times (x,\Bbb
C)@>{\mu_x}>>(M,\mathcal O_{M})
\end{CD}.
$$
Since $(\iota,\widehat{x})$ is an isomorphism, we get the
commutativity of (\ref{iota(G_x)}).

Denote by $\mathcal{I}_{\iota (G_x)}$ the sheaf of ideals in
$\mathcal O_{G}$ corresponding to the subsupermanifold\linebreak
$(\iota_1(G_x),\mathcal O_{\iota_1(G_x)})\subset (G,\mathcal
O_{G})$. Using (\ref{iota(G_x)}), we have:
$$
\phi^i_U(x_s,\xi_t)+\mathcal{I}_{\iota
(G_x)}=\mu_x(y_i)+\mathcal{I}_{\iota (G_x)} =
y_i(x)+\mathcal{I}_{\iota (G_x)}=0+\mathcal{I}_{\iota
(G_x)}=\mathcal{I}_{\iota (G_x)}.
$$
Therefore $\phi^i_U(x_s,\xi_t)\in \mathcal{I}_{\iota (G_x)}$.
Similarly, $\psi^i_U(x_s,\xi_t)\in \mathcal{I}_{\iota (G_x)}$. It
follows that $\mathcal{I}\subset \mathcal{I}_{\iota (G_x)}
=\iota_2(\mathcal{I})$, and the equality $\iota^2=\id$ implies that
$\mathcal{I}= \mathcal{I}_{\iota (G_x)}$.$\Box$

\bigskip

We have seen that for a transitive action of $(G,\mathcal O_{G})$ on
$(M,\mathcal O_{M})$ the ringed space $(G_x,\mathcal O_{G_x})$ is a
supermanifold, which is in fact a Lie subsupergroup of $(G,\mathcal
O_{G})$. This Lie subsupergroup will be called
 {\it the stationary
 Lie subsupergroup of $x$}. In the last section we will prove that any homogeneous supermanifold is
isomorphic to a coset space of a Lie supergroup with the structure
of a supermanifold introduced above.

\bigskip

\bigskip

\begin{center}
{\bf 4. The structure of a homogeneous space}
\end{center}

\medskip

Our goal here is the following theorem.

\medskip

\t\label{G/G_m_i_(M,O)} {\sl Let $(M,\mathcal O_M)$ be a
$(G,\mathcal O_G)$-homogeneous supermanifold and let $(H,\mathcal
O_H)$ be a stationary Lie subsupergroup  of a point $x\in M$. Then
there is an isomorphism
$$
\beta: (G/H, \mathcal O_{G/H})\to (M,\mathcal O_M),
$$
which is $(G,\mathcal O_G)$-equivariant in the sense that the
following diagram is commutative:
\begin{equation}
\label{ekvivariantnost'}
\begin{CD}
(G,\mathcal O_G)\times (G/H, \mathcal O_{G/H})@>{\id\times \beta}>>
(G,\mathcal O_G)\times (M,\mathcal O_M)
\\
@V{\alpha}VV @VV{\mu}V\\
(G/H, \mathcal O_{G/H})@>{\beta}>>  (M,\mathcal O_M)
\end{CD}.
\end{equation}
 }
\medskip

\noindent{\it Proof}. Let $\beta_1:G/H\to M$, $gH\mapsto gx$, be the
natural homeomorphism. Then the diagram
\begin{equation}
\label{digramma_M_isom_G/G_m}
\begin{CD}
(G,\mathcal O_G)@= (G,\mathcal O_G)
\\
@V{p=(p_1,p_2)}VV @VV{\mu_x=((\mu_x)_1,(\mu_x)_2),}V\\
(G/H, \mathcal O_{G/H})@>{\beta_1}>>  (M,\mathcal O_M)
\end{CD}
\end{equation}
is commutative. We will now construct an isomorphism of sheaves
$\beta_2:\mathcal O_M\to (\beta_1)_*\mathcal O_{G/H}$, such that
$\beta = (\beta _1, \beta _2)$ is the required isomorphism.


As in Theorem \ref{triv rassl}, consider the product $(S,\mathcal
O_S)\times(H,\mathcal O_{H})$. By the axioms of action, the
composition of morphisms
$$
\begin{CD}
(S,\mathcal O_S)\times(H,\mathcal O_{H})@>{\nu|_{S\times H}}>>
(U,\mathcal O_G)@>{\mu_x}>> (M,\mathcal O_M).
\end{CD}
$$
can be written as $\mu_x\circ \nu|_{S\times H}= \mu\circ (\id \times
\mu_x|_H)=\mu \circ (\id \times \widehat{x})$. The last equality
follows from the fact that $\mu_x|_{H}=\widehat{x}$. As a result, we
get
\begin{equation}
\label{zvezda}
\mu_x \circ \nu|_{S\times H}= \mu_x\circ \pr^{S\times
H}_{S}.
\end{equation}
 The
differential $(\d\!\mu_x)_{e}$ is surjective, because the action
$\mu$ is transitive. The differential $(\d\!\nu|_{S\times
H})_{(e,e)}$ is nondegenerate, see Theorem \ref{triv rassl}. It is
easy to see that $\dim(S,\mathcal O_S)= \dim(M,\mathcal O_M)$, and
so we get that $(\d\!(\mu_x)|_{S})_{(e)}$ is nondegenerate. By
Theorem \ref{ob obratnoj funkcii}, we can assume that $\mu|_S$ is an
isomorphism of $(S,\mathcal O_S)$ onto some superdomain $(V,\mathcal
O_M)\subset(M,\mathcal O_M)$. By (\ref{zvezda}) we get
$(\mu_x)_2(\mathcal O_M|_V)=((\nu|_{S\times H})^{-1})_2(\mathcal
O_S)$. On the other hand, by the definition of $p$ we have
$p_2(\mathcal O_{G/H}|_{p_1(S\times H)})=((\nu|_{S\times
H})^{-1})_2(\mathcal O_S)$ and therefore $p_2(\mathcal
O_{G/H}|_{p_1(S\times H)})=(\mu_x)_2(\mathcal O_M|_V)$.

Let $\widetilde{l}_g=\alpha \circ (\widehat{g}\times \id)$. The
following diagram is commutative by Lemma \ref{sdvig_deystvija}:
\begin{equation}
\label{digrammap,l_g}
\begin{CD}
(p_1(S\times H),\mathcal O_{G/H})@>{\widetilde{l}_g}>> (p_1(gS\times
H),\mathcal O_{G/H})
\\
@A{p}AA @A{p}AA\\
(U,\mathcal O_{G})@>{l_g}>>
(gU,\mathcal O_{G})\\
@V{\mu_x}VV @V{\mu_x}VV\\
(V,\mathcal O_M)@>{\overline{l}_g}>>  (gV,\mathcal O_M)
\end{CD}.
\end{equation}
By the commutativity of (\ref{digrammap,l_g}), we get
$$
(\mu_x)_2(\mathcal O_M|_{gV})=p_2(\mathcal O_{G/H}|_{p_1(gU)}),\quad
\text{for all}\  g\in G.
$$
Therefore, the sheaves $(\mu_x)_2(\mathcal O_{M})$ and $p_2(\mathcal
O_{G/H})$ coincide locally, and it follows that $p_2(\mathcal
O_{G/H})=(\mu_x)_2(\mathcal O_{M})$. Define a morphism
$\beta_2:\mathcal O_{M}\to (\beta_1)_*(\mathcal O_{G/H})$ by
$$
\beta_2(f)=p_2^{-1}\circ (\mu_x)_2(f).
$$
Obviously, the morphism $\beta=(\beta_1,\beta_2)$ is an isomorphism
and the diagram (\ref{digramma_M_isom_G/G_m}) is commutative.

It remains to prove that diagram (\ref{ekvivariantnost'}) is also
commutative.  By the definition of $\beta$, we get locally
$\beta=(\beta_1,\beta_2)=\mu_x\circ i_g$, where $i_g$ is defined
after Lemma \ref{f-(H,O_H)-invariantna2}. 
 Hence $\mu_x\circ
i_g\circ p=\mu_x$ and the axioms of action yield $\mu\circ
(\id\times \mu_x)=\mu_x\circ \nu$. Thus we get locally:
$$
\mu\circ (\id\times \beta)= \mu\circ (\id\times \mu_x\circ i_g)=
\mu\circ (\id\times \mu_x)\circ (\id\times i_g)=\mu_x\circ \nu \circ
(\id\times i_g).
$$
$$
\beta \circ\alpha=\mu_x\circ i_g\circ p \circ\nu\circ (\id\times
i_g)= \mu_x \circ\nu\circ (\id\times i_g).
$$
This implies the commutativity of the diagram
(\ref{ekvivariantnost'}).$\Box$

\bigskip

\emph{The author would like to thank her supervisor Prof. Arkady
Onishchik for many helpful remarks on the paper, Prof. Dmitry
Akhiezer for interesting discussions and Prof. Dr. Dr. h.c. mult.
Alan T. Huckleberry for the hospitality at Ruhr-Universit\"{a}t
Bochum.}

\smallskip

\textsc{Tver State University, Zhelyabova 33, 170 000 Tver, Russia}

\emph{E-mail address:} \verb"VishnyakovaE@googlemail.com"


\medskip


\begin{thebibliography}{99}

\bibitem{Kostant} {\it Kostant B.} Graded Manifolds, graded Lie
theory, and prequantization. Lecture Notes in Mathematics 570.
Berlin e.a.: Springer-Verlag,  1977. P. 177-306.

\bibitem{ley} {\it Leites D.A.} Introduction to the theory of supermanifolds.
Russian Math. Surveys 35 (1980), 1-64.

\bibitem{man} {\it Manin Yu. I.} Gauge field theory and complex
geometry. Grundlehren der Mathematischen Wissenschaften, 289,
Springer-Verlag, Berlin, 1997.


\bibitem{onipi} {\it Onishchik A.L.} Flag supermanifolds, their
automorphisms und deformations. The Sophus Lie Memorial conference
(Oslo, 1992), 289-302, Scand. Univ. Press, Oslo, 1994.




\end{thebibliography}
\end{document}